\documentclass[10.5pt, a4paper,  reqno, english]{amsart}
\usepackage[all]{xy}
\usepackage{amssymb,amsmath,amsthm}
\numberwithin{equation}{section}

\newtheorem{theorem}{Theorem}
\newtheorem{proposition}{Proposition}

\theoremstyle{definition}

\newcommand{\ben}{\begin{enumerate}}
\newcommand{\een}{\end{enumerate}}
\newcommand{\eit}{\begin{itemize}}
\newcommand{\beq}{\begin{equation}}
\newcommand{\eeq}{\end{equation}}

\renewcommand{\leq}{\leqslant}
\renewcommand{\geq}{\geqslant}

\begin{document}

\title{An alternative to Vaughan's identity}

\author{ Andrew Granville}

 \address{D{\'e}partment  de Math{\'e}matiques et Statistique,   Universit{\'e} de Montr{\'e}al, CP 6128 succ Centre-Ville, Montr{\'e}al, QC  H3C 3J7, Canada; and  Department of Mathematics, University College London, Gower Street, London WC1E 6BT, England.}
   \email{andrew@dms.umontreal.ca}  
   \thanks{The author has received funding from the
European Research Council  grant agreement n$^{\text{o}}$ 670239, and from NSERC Canada under the CRC program.}
\begin{abstract}
We exhibit an identity that plays the same role as Vaughan's identity but is arguably simpler.
\end{abstract}

\maketitle

\newcommand{\cbar}{\overline{\chi}}
\newcommand{\pbar}{\overline{\psi}}
\newcommand{\sumstar}{\sideset\and ^* \to \sum}

\section{Introduction}

Let $1_y$ denote the characteristic function of the integers free of prime factors $\leq y$,
The idea is to work with the   identity
\begin{equation}  \label{eq: TrivId}
\Lambda(n) = \log n - \sum_{\substack{\ell m=n \\ \ell,m>1 }} \Lambda(\ell) ,
\end{equation}
summing it up   over  integers $n\leq x$ for which $1_y(n)=1$, where we might select $y:=\exp(\sqrt{\log x})$ or larger. 
In this case the second sum can be written a sum of terms $1_y(\ell) \Lambda(\ell) \cdot 1_y(m)$ which have the bilinear structure that is used in ``Type II sums''. We will see the identity in action in two key results in analytic number theory:

\section{The Bombieri-Vinogradov Theorem}

\begin{theorem}  [The Bombieri-Vinogradov Theorem] \label{thm: 46.3}
If $x^{1/2}/(\log x)^B\leq Q\leq x^{1/2}$ then
  \begin{equation}  \label{eq: 46.6}
\sum_{q\leq Q} \max_{(a,q)=1} \ \left| \pi(x;q,a) - \frac{\pi(x)}{\phi(q)}  \right|  \ll Q x^{1/2}  (\log\log x)^{1/2} .
\end{equation}
\end{theorem}

This is a little stronger than the results in the literature (for example Davenport \cite{Da} has the 
$(\log\log x)^{1/2}$ replaced by $(\log x)^{5}$). The reason for this improvement is the simplicity of our identity, and some slight strengthening of the auxiliary results used in the proof.

\begin{proof} Let $y=x^{1/\log\log x}$ . We will instead prove the following result, in which the $\psi$ function replaces $\pi$, and deduce \eqref{eq: 46.6} by partial summation:
\begin{equation}  \label{eq: 46.5}
\sum_{q\leq Q} \max_{(a,q)=1} \ \left| \psi(x;q,a) - \frac{\psi(x)}{\phi(q)}  \right|  \ll Q x^{1/2}  (\log x) (\log\log x)^{1/2} .
\end{equation}

Using \eqref{eq: TrivId} for integers $n$ with $1_y(n)=1$, the quantity on the left-hand side of \eqref{eq: 46.5} is $\leq S_I+S_{II}+E$ where
\[
S_I=\sum_{q\leq Q} \max_{(a,q)=1} \ \left| \sum_{\substack{ n\leq x \\   n\equiv a \pmod q  }} 1_y(n) \log n - \frac {1}{\phi(q)} \sum_{\substack{ n\leq x \\   (n,q)=1  }} 1_y(n) \log n  \right|   
\]
which is $\ll xu^{-u+o(u)}\ll_A \frac{x}{(\log x)^A}$ by the small sieve, where $x/Q=y^u$; 
and $E$ is the contribution of the powers of primes $\leq y$, which contribute $\leq \pi(y)\log x$ to each sum and therefore $\leq Q \pi(y)\log x\ll_A \frac{x}{(\log x)^A}$ in total. Most interesting is
\[
S_{II}=\sum_{q\leq Q} \ \max_{a:\ (a,q)=1} \ \left| \sum_{\substack{ n\equiv a \pmod q  }} f(n) - \frac {1}{\phi(q)} \sum_{\substack{ (n,q)=1  }} f(n)\right| 
\]
where $f(n)= \sum_{\substack{\ell m=n, \ell,m>y }}  \Lambda(\ell)1_y(\ell) \cdot 1_y(m)$. Its bilinearity means that this is a Type II sum, and we can employ the following general result.\footnote{This can be obtained by taking the ideas in proving (5) of chapter 28 of \cite{Da},  along with the method of proof of Theorem 9.16 of \cite{FI}; in any case it is only a minor improvement on either of these results. For full details see chapter 51 of \cite{Gra}.}

\begin{theorem} \label{thm: 46.1} For each integer $n\leq x$ we define
\[
f(n):= \sum_{\substack{\ell m=n }}  \alpha_\ell\beta_m
\]
where $\{ \alpha_\ell\}$ and $\{ \beta_m\}$ are sequences of complex numbers, for which  
\begin{itemize} 
\item The  $\{ \alpha_\ell\}$ satisfy the Siegel-Walfisz criterion;
\item The  $\{ \alpha_\ell\}$  are only supported in the range $L_0 \leq  \ell  \leq x/y$ ;
\item $\sum_{\ell \leq L} |\alpha_\ell |^2 \leq aL$ and $\sum_{m \leq M} |\beta_m|^2 \leq bM$ for all $L,M\leq x$.
\end{itemize}
For any $B>0$ we have
\begin{equation}  \label{eq: 46.4}
\sum_{q\leq Q} \ \max_{a:\ (a,q)=1} \ \left| \sum_{\substack{ n\equiv a \pmod q  }} f(n) - \frac {1}{\phi(q)} \sum_{\substack{ (n,q)=1  }} f(n)\right| \ll (ab)^{1/2} Q x^{1/2}  \log x,
\end{equation}
where $Q=x^{1/2}/(\log x)^B$,  with $x/y\leq  \frac{Q^2}{(\log x)^2}$ and $L_0\geq y,\exp( (\log x)^\epsilon)$.
\end{theorem}

We deduce that 
$S_{II} \ll Q x^{1/2}  (\log x)^{5/4} $ by Theorem \ref{thm: 46.1} since $a\ll \log x$ and $b\ll \frac 1{\log y}=  \frac{\log\log x}{\log x}$.
  \end{proof}

\section{A general bound for a sum over primes}
 
\begin{proposition} \label{prop: 80.General} For any given function $F(.)$ and $y\leq x$ we have 
\[
\left|\sum_{ \substack{ n\leq x \\  p(n)>y}}\Lambda(n)F(n)\right| \ll S_I \log x + (S_{II} \, x (\log x)^5)^{1/2}
\]
where $S_I$ is the Type I sum given by 
\[
S_I:= \max_{t\leq x} \left| \sum_{ \substack{   n\leq t \\  p(n)>y}}   F(n) \right|  \leq   \sum_{\substack{   d\geq 1 \\  P(d)\leq y}} \left| \sum_{   m\leq t/d }   F(d  m) \right| ,
\]
and $S_{II}$ is the Type II sum given by 
\[
S_{II}:= \max_{\substack{ y<L\leq x/y \\ y< m\leq 2x/L}}   \sum_{m/2<n\leq 2m}  \left| \sum_{\substack{ L<\ell\leq 2L \\ \ell \leq  \frac xm ,\frac xn }} F(\ell m)\overline{F(\ell n)} \right|
\]
\end{proposition}

This simplifies, and slightly improves chapter 24 of \cite{Da}, which is what is used there to bound exponential sums over primes.

\begin{proof} We again use  \eqref{eq: TrivId} so that 
\[
\sum_{ \substack{ n\leq x \\  p(n)>y}} \Lambda(n)F(n) =
\sum_{ \substack{ n\leq x \\  p(n)>y}} F(n)    \log n -  \sum_{\substack{ \ell  \\  p(\ell)>y}} \Lambda(\ell)   \sum_{\substack{ m\leq x/\ell \\  p(m)>y}}  F(\ell m).
\]
where $p(n)$ denotes the smallest prime factor of $n$. Now
\[
\sum_{ \substack{ n\leq x \\  p(n)>y}} F(n)    \log n =\sum_{ \substack{ n\leq x \\  p(n)>y}} F(n) \int_1^n \frac {dt}t =  \int_1^x  \sum_{ \substack{ t\leq n\leq x \\  p(n)>y}}    F(n) \frac {dt}t \leq
2\log x \cdot \max_{t\leq x} \left| \sum_{ \substack{   n\leq t \\  p(n)>y}}   F(n) \right|.
\]
Moreover for $P=\prod_{p\leq y}p$,
\[
 \sum_{  \substack{   n\leq t \\  p(n)>y}}   F(n)  =  \sum_{   n\leq t }   F(n)    \sum_{d|P, d|n} \mu(d)
 =  \sum_{d|P } \mu(d) \sum_{   m\leq t/d }   F(d m) .
\]
 
For the second sum we first split the sums into dyadic intervals ($L<\ell\leq 2L, M<m\leq 2M$) and then Cauchy, so that the square of each subsum is
\begin{align*}
& \leq \sum_{\ell : p(\ell)>y} \Lambda(\ell)^2 \cdot \sum_\ell \left|  \sum_{\substack{ m\leq x/\ell \\  p(m)>y}}    F(\ell m)  \right|^2
\ll L \log L \sum_{\substack{M<m,n\leq 2M  \\  p(m),p(n)>y}}  \sum_{\ell \leq \frac x{\max\{ m,n\} }} F(\ell m)\overline{F(\ell n)} \\ 
& \ll  x\log x \cdot   \max_{M<m\leq 2M}  \sum_{\substack{m/2<n\leq 2m\\  p(n)>y}}  \left| \sum_{\substack{ L<\ell\leq 2L \\ \ell \leq  \frac xm ,\frac xn }} F(\ell m)\overline{F(\ell n)} \right| \end{align*}
since $m,n\in (M,2M]$, and the result follows.
\end{proof}


\section{Genesis}

The idea for using \eqref{eq: TrivId} germinated from reading the proof of the Bombieri-Vinogradov Theorem (Theorem 9.18)  in \cite{FI}, in which they used Ramar\'e's identity, that if $\sqrt{x}<n\leq x$ and $n$ is squarfree then
\[
  1_{\mathbb P}(n) =  1 - \sum_{\substack{pm=n \\ p \text{ prime} \leq \sqrt{x}}} \frac 1{\omega_{\sqrt{x}}(m)}
\]
where $1_{\mathbb P}$ is the characteristic function for the primes, and  $\omega_{z}(m)=1+\sum_{p|m,\ p\leq z} 1$.
They also had to sum this over all integers free of prime factors $>y$.

\bibliographystyle{plain}

\end{document}